\documentclass[12pt,oneside,reqno]{amsart}
 \usepackage{amscd,amsmath,amsthm,amssymb}
 \def\NZQ{\Bbb}               % the font for N,Z,Q,R,C
 \def\NN{{\NZQ N}}

 \def\ab{{\bold a}}
 
 \def\xb{{\bold x}}

 \def\opn#1#2{\def#1{\operatorname{#2}}} % to make operators
 \opn\chara{char} \opn\length{\ell} \opn\pd{pd} \opn\rk{rk}
 \opn\projdim{proj\,dim} \opn\injdim{inj\,dim} \opn\rank{rank}
 \opn\depth{depth} \opn\grade{grade} \opn\height{height}
 \opn\embdim{emb\,dim} \opn\codim{codim}
 
 \opn\Tr{Tr} \opn\bigrank{big\,rank}
 \opn\superheight{superheight}\opn\lcm{lcm}
 \opn\trdeg{tr\,deg}%\emph{
 \opn\reg{reg} \opn\lreg{lreg} \opn\ini{in} \opn\lpd{lpd}
 \opn\size{size} \opn\sdepth{sdepth}
 \opn\link{link}\opn\fdepth{fdepth}\opn\lex{lex}
 \opn\div{div} \opn\Div{Div} \opn\cl{cl} \opn\Cl{Cl}
 \opn\Spec{Spec} \opn\Supp{Supp} \opn\supp{supp} \opn\Sing{Sing}
 \opn\Ass{Ass} \opn\Min{Min}\opn\Mon{Mon}
 \opn\Ann{Ann} \opn\Rad{Rad} \opn\Soc{Soc}
 \opn\Im{Im} \opn\Ker{Ker} \opn\Coker{Coker} \opn\Am{Am}
 \opn\Hom{Hom} \opn\Tor{Tor} \opn\Ext{Ext} \opn\End{End}
 \opn\Aut{Aut} \opn\id{id}
 
 \opn\nat{nat}
 \opn\pff{pf}%   \pf exists already
 \opn\Pf{Pf} \opn\GL{GL} \opn\SL{SL} \opn\mod{mod} \opn\ord{ord}
 \opn\Gin{Gin} \opn\Hilb{Hilb}\opn\sort{sort}
 \opn\Tot{Tot}
 \opn\aff{aff} \opn
 \con{conv} \opn\relint{relint} \opn\st{st}
 \opn\lk{lk} \opn\cn{cn} \opn\core{core} \opn\vol{vol}
 \opn\link{link} \opn\star{star}\opn\lex{lex}\opn\set{set}
 \opn\gr{gr}
 
 \def\pot#1#2{#1[\kern-0.28ex[#2]\kern-0.28ex]}

 \opn\dirlim{\underrightarrow{\lim}}
 \opn\inivlim{\underleftarrow{\lim}}
 \let\union=\cup

 \def\Implies{\ifmmode\Longrightarrow \else
         \unskip${}\Longrightarrow{}$\ignorespaces\fi}
 \def\implies{\ifmmode\Rightarrow \else
         \unskip${}\Rightarrow{}$\ignorespaces\fi}
 \def\iff{\ifmmode\Longleftrightarrow \else
         \unskip${}\Longleftrightarrow{}$\ignorespaces\fi}

 \let\:=\colon
 \newtheorem{Theorem}{Theorem}[section]

 \newtheorem{Proposition}[Theorem]{Proposition}

 \newtheorem{Example}[Theorem]{Example}

 \let\epsilon\varepsilon
 \let\kappa=\varkappa
 \def\qed{\ifhmode\textqed\fi
       \ifmmode\ifinner\quad\qedsymbol\else\dispqed\fi\fi}
 \def\textqed{\unskip\nobreak\penalty50
        \hskip2em\hbox{}\nobreak\hfil\qedsymbol
        \parfillskip=0pt \finalhyphendemerits=0}
 \def\dispqed{\rlap{\qquad\qedsymbol}}
 
 \opn\dis{dis}
 \def\pnt{{\raise0.5mm\hbox{\large\bf.}}}
 
 \opn\Lex{Lex}

\begin{document}

\title{Polymatroidal property of generalized mixed product ideals}

\author {Monica La Barbiera and Roya Moghimipor}

\address{Monica La Barbiera, Department of Mathematics and Informatics, University of Messina, Viale Ferdinando
Stagno d’Alcontres 31, 98166 Messina, Italy} \email {monicalb@unime.it}

\address{Roya Moghimipor, Department of Mathematics, Safadasht Branch, Islamic Azad University, Tehran, Iran} \email{roya\_moghimipour@yahoo.com}

\begin{abstract}
Let $L$ be the generalized mixed product ideal induced by a monomial ideal $I$.
In this paper, we study the polymatroidal property of generalized mixed product ideals.
Furthermore, some algebraic invariants of $L$ are computed.
\end{abstract}


\maketitle

\section{Introduction}
\label{one}
Restuccia and Villarreal \cite{RV} introduced mixed product ideals, which form  a particular class of squarefree monomial ideals.
They gave a complete classification of normal mixed product ideals, as well as applications
in graph theory.

Rinaldo \cite{Ri} and Ionescu and Rinaldo \cite{IR} studied algebraic and homological properties of this class of ideals, and Hoa and Tam \cite{HT} computed the regularity and some other algebraic invariants of mixed products of arbitrary graded ideals.
In \cite{Ri}, Rinaldo classified the ideals of mixed products that are sequentially Cohen-Macaulay. In \cite{LR} the first author together with Restuccia investigated monomial ideals of mixed products generated by a s-sequence in order to compute invariants of their symmetric algebra.

Let $K$ be a field and $K[x_1,\dots,x_n]$ the polynomial
ring in $n$ variables over $K$ with each $x_i$ of degree $1$.
Let $I\subset S$ be a monomial ideal and $G(I)$ its unique
minimal monomial generators.

Let $K$ be a field and  $S = K[x_1, \ldots,  x_n, y_1, \ldots, y_m]$  be the polynomial ring over $K$ in the variables $x_i$ and $y_j$.

Mixed product ideals are of the form $(I_qJ_r + I_pJ_s)S$, where for integers $a$ and $b$, the ideal  $I_a$  (resp.\ $J_b$) is the ideal generated by all squarefree monomials of degree $a$ in the polynomial ring $K[x_1,\ldots, x_n]$ (resp.\ of degree $b$ in the polynomial ring $K[y_1, \ldots, y_m]$), and where $0 < p < q \leq n$, $0 < r < s \leq m$.
Thus the ideal $L=(I_qJ_r + I_pJ_s)S$ is obtained from the monomial ideal $I=(x^qy^r,x^py^s)$ by replacing $x^q$ by $I_q$, $x^p$ by $I_p$, $y^r$ by $J_r$ and $y^s$ by $J_s$.

Together with  Herzog and  Yassemi \cite{HMY} the second author introduced the generalized mixed product ideals, which  are a far reaching generalization of  the mixed product ideals introduced by Restuccia and Villarreal, and also generalizes the expansion construction by Bayati and Herzog \cite{BH}.
A great deal of knowledge on the generalized mixed product ideals is accumulated in several papers \cite{HMY, Mo, MN, MB}.

For this new construction we choose for each $i$ a set of new variables $x_{i1},x_{i2},\ldots,x_{im_i}$ and replace each of the factor $x_i^{a_i}$ in each minimal generator $x_1^{a_1}x_2^{a_2}\cdots x_n^{a_n}$ of the monomial ideal $I$ by a monomial ideal in $T_i=K[x_{i1},x_{i2},\ldots,x_{im_i}]$ generated in degree $a_i$.

We computed in \cite{HMY} the minimal graded free resolution of generalized mixed product ideals and showed that a generalized mixed product ideal $L$ induced by $I$ has the same regularity as $I$, provided the ideals which replace the pure powers $x_i^{a_i}$ all have a linear resolution. As a consequence we obtained the result that under the above assumptions, $L$ has a linear resolution if and only if $I$ has a linear resolution. We also proved that the projective dimension of $L$ can be expressed in terms of the multi-graded shifts in the resolution of $I$ and the projective dimension of the ideals which replace the
pure powers.

In \cite{M} the second author together with Tehranian introduced the generalized mixed polymatroidal ideals. The class of generalized mixed polymatroidal ideals is a particular class of generalized mixed product ideals which for each $i$ we replace each factor $x_i^{a_i}$ in each minimal generator $x_1^{a_1}x_2^{a_2}\cdots x_n^{a_n}$ of $I$ by a polymatroidal ideal in $T_i$ generated in degree $a_i$. Also we computed powers of generalized mixed product ideals and showed that
$L^k$ is again generalized mixed product ideal for all $k$ and $L^k$ induced by $I^k$, and we obtained the result that $L^k$ has a linear resolution if and
only if $I^k$ has a linear resolution for all $k$, provided the ideals which replace the pure powers $x_i^{a_i}$ all have a linear resolution.

For a given positive integer $k$, we denote the $k$-th bracket power
of $I$, to be the ideal $I^{[k]}$, generated by all monomials $u^{k}$, where
$u\in I$ is a monomial.
In \cite{Mo}, we characterized the Cohen-Macaulay ideals of generalized mixed products and, in \cite{MB} algebraic properties of $L^{[k]}$ are studied.

The present paper is organized as follows. In Section 2 we study the polymatroidal property of generalized mixed product ideals. We present several cases for  which  a generalized mixed product ideal $L$ induced by the monomial ideal $I$ is polymatroidal.
In Theorem \ref{normality} we study the normality of $L$ induced by a monomial ideal in $K[x_1,x_2]$,
where the ideals who substitute the generators of $I$ are squarefree Veronese.

In Section \ref{three} we study the linear quotients of generalized mixed product ideals, and some invariants of $L$ are computed.

\section{Polymatroidal ideals}
\label{two}
In this section we want to study the polymatroidal property of the generalized mixed product ideals.
Let $S=K[x_{1},\dots,x_{n}]$ be the polynomial ring over a field $K$ in the variables $x_1,\ldots,x_n$ with the maximal ideal $\mathfrak{m}=(x_1,\dots,x_n)$, and let $I\subset S$ be a monomial ideal with $I\neq S$ whose minimal set of generators  is $G(I)=\{\xb^{\ab_1},\ldots, \xb^{\ab_m}\}$. Here $\xb^{\ab_{j}}=x_1^{\ab_{j}(1)}x_2^{\ab_{j}(2)}\cdots x_n^{\ab_{j}(n)}$ for $\ab_{j}=(\ab_{j}(1),\ldots,\ab_{j}(n))\in\NN^n$.

Next we consider the polynomial ring $T$ over $K$ in the variables
\[
x_{11},\ldots,x_{1m_1},x_{21},\ldots,x_{2m_2},\ldots,x_{n1},\ldots,x_{nm_n}.
\]
The class of polymatroidal ideals is one of the rare classes of monomial ideals
with the property that all powers of an ideal in this class have a linear resolution.

Recall that a monomial ideal is called polymatroidal, if its monomial generators correspond to the bases of a discrete
polymatroid, see \cite{HH}. Since the set of bases of a discrete polymatroid is characterized
by the so-called exchange property, it follows that a polymatroidal ideal may as
well be characterized as follows: let $I\subset S$ be a monomial ideal generated in a single degree and $G(I)$ its unique
minimal monomial generators.
Thus $I$ is said to be {\em polymatroidal}, if for any two elements $u, v \in G(I)$ such that $\deg_{x_i}(u)> \deg_{x_i}(v)$
there exists an index $j$ with $\deg_{x_j}(u)< \deg_{x_j}(v)$ such that $x_{j}(u/x_{i})\in I$.

In \cite{HMY} we introduced the generalized mixed product ideals.
For $i=1,\ldots,n$ and $j=1,\ldots,m$ let $L_{i, \ab_j(i)}$ be a monomial ideal in the variables $x_{i1},x_{i2},\ldots,x_{im_i}$ such that
\begin{eqnarray}
\label{inclusion}
L_{i, \ab_j(i)}\subset L_{i, \ab_k(i)} \quad \text{whenever} \quad \ab_j(i)\geq \ab_k(i).
\end{eqnarray}

Given these ideals we define for $j=1,\ldots,m$ the monomial ideals
\begin{eqnarray}
\label{lj}
L_j=\prod_{i=1}^nL_{i, \ab_j(i)} \subset T,
\end{eqnarray}
and set $L=\sum_{j=1}^mL_j$.
The ideal $L$ is called a {\em generalized mixed product ideal} induced by $I$.

\begin{Example}
{\em
Consider the mixed product ideals introduced by Restuccia and Villarreal \cite{RV}.
A mixed product ideal is a monomial ideal of the the form
\[
L=I_qJ_r + I_sJ_t,
 \]
where for integers $a$ and $b$, the ideal  $I_a$  (resp.\ $J_b$) is the ideal generated by all squarefree monomials of degree $a$ in the polynomial ring $K[x_1,\ldots, x_n]$ (resp.\ of degree $b$ in the polynomial ring $K[y_1, \ldots, y_m]$), and where $0 \leq s < q \leq n$, $0 \leq r < t \leq m$.
Ideals of this type are called squarefree Veronese ideals.
}
\end{Example}

A generalized mixed product ideal depends not only on $I$ but also on the family $L_{ij}$.
Thus we write $L(I;\{L_{ij}\})$ for the generalized mixed product ideal induced by a monomial ideal $I$.

Now we assume that the generalized mixed product ideal $L$ induced by a monomial ideal $I$, where the ideals $L_{ij}$ are squarefree Veronese ideals in $K[x_{i1},x_{i2},\ldots,x_{im_i}]$ of degree $a_{j}$.

\begin{Theorem}
\label{qr}
Let $L=I_{q}J_{r}\subset K[x_1,\ldots, x_n,y_1, \ldots, y_m]$ be the generalized mixed product ideal, where for integers $a$ and $b$, the ideal $I_a$ (resp.\ $J_b$) is the ideal generated by all squarefree monomials of degree $a$ in the polynomial ring $K[x_1,\ldots, x_n]$
(resp.\ of degree $b$ in the polynomial ring $K[y_1, \ldots, y_m]$).
Then $L$ is matroidal if and only if one of the following conditions are satisﬁed:

{\em (1)} $q=0$;

{\em (2)} $r=0$;

{\em (3)} $q,r> 0$.
\end{Theorem}

\begin{proof}
The cases $q=0$ or $r=0$ being clear, we can suppose that $q>0$ and $r>0$.
Since the squarefree Veronese ideal is matroidal, and then by applying \cite[Theorem 12.6.3]{HH} we get that $I_{q}J_{r}$ is polymatroidal, as desired.
\end{proof}

\begin{Theorem}
\label{polymatroidal}
Let $L=I_{q}J_{r}+I_{q+1}J_{r-1} \subset K[x_1,\ldots, x_n,y_1, \ldots, y_m]$ be the generalized mixed product ideal, where for integers $a$ and $b$, the ideal $I_a$  (resp.\ $J_b$) is the ideal generated by all squarefree monomials of degree $a$ in the polynomial ring $K[x_1,\ldots, x_n]$
(resp.\ of degree $b$ in the polynomial ring $K[y_1, \ldots, y_m]$), and where $q \geq0$, $r\geq 1$. Then $L$ is matroidal.
\end{Theorem}

\begin{proof}
Since
\[
G(L)=G(I_{q}J_{r}+I_{q+1}J_{r-1})=G(I_{q}J_{r})\union G(I_{q+1}J_{r-1}),
\]
we consider $u,v$ be two monomials belonging to $G(L)$.
We distinguish several cases:

$(\.{i})$ If $u,v \in G(I_{q}J_{r})$, then $L$ is matroidal by Theorem \ref{qr}.

$(\.{i}\.{i})$   If $u,v\in G(I_{q+1}J_{r-1})$, thus Theorem \ref{qr} yields $L$ is matroidal.

$(\.{i}\.{i}\.{i})$ Let
\[
u=x_{i_{1}}\dots x_{i_{q}} y_{j_{1}}\dots y_{j_{r}}\in G(I_{q}J_{r}),
\]
\[
v= x_{i_{1}}\dots x_{i_{q+1}}y_{j_{1}}\dots y_{j_{r-1}}\in G(I_{q+1}J_{r-1})
\]
such that $i_{1} < i_{2} < \dots i_{q}<i_{q+1}$ and $j_{1}< j_{2}<\dots <j_{r-1}< j_{r}$.

If $\deg_{y_{j_{p}}}(u)>\deg_{y_{j_{p}}}(v)$, then there exists an integer $q+1$ such that $\deg_{x_{i_{q+1}}}(u)< \deg _{x_{i_{q+1}}}(v)$ and
$x_{i_{q+1}}(u/y_{j_{p}})\in G(L)$.

If $\deg_{x_{i_{l}}}(v)>\deg_{x_{i_{l}}}(u)$, then there exists an integer $r$ such that $\deg_{y_{i_{r}}}(v)< \deg _{y_{i_{r}}}(u)$ and
$y_{j_{r}}(u/x_{i_{l}})\in G(L)$. Then the desired conclusion follows.

\end{proof}

\begin{Theorem}
Let $L=I_{1}+J_{1}$ be the generalized mixed product ideal, where the ideal $I_1$  (resp.\ $J_1$) is the ideal generated by all squarefree monomials of degree 1 in the polynomial ring $K[x_1,\ldots, x_n]$ (resp.\ of degree 1 in the polynomial ring $K[y_1, \ldots, y_m]$).
Then $L$ is matroidal.
\end{Theorem}

\begin{proof}
(a) Let $L=I_{1}+J_{1} \subset K[x_1,\dots,x_n, y_1,\dots,y_m]$ be a generalized mixed product ideal induced by the ideal $I=(x,y)\subset K[x,y]$.
Therefore,
\[
L=(x_1,\dots,x_n,y_1,\dots,y_m)
\]
is a maximal ideal, and hence $L$ is matroidal.
\end{proof}

\begin{Example}
\label{pathideal}
{\em
As mentioned in \cite{RV} mixed product ideals also appear as generalized graph ideals (called path ideals by Conca and De Negri \cite{CN}) of complete bipartite graphs. Let $G$ a finite simple graph with vertices $x_1,\ldots,x_n$. A {\em path} of length $t$ in $G$ is sequence $x_{i_1},\ldots,x_{i_t}$ of pairwise distinct vertices such that $\{x_{i_k},x_{i_{k+1}}\}$ is an edge of  $G$. Then the  {\em path ideal} $I_t(G)$ is the ideal generated by all monomials $x_{i_1}\cdots x_{i_t}$ such that  $x_{i_1},\ldots,x_{i_t}$ is a path of length $t$.

Now let $G$ be a complete bi-partite graph with vertex set $V=V_1\union V_2$, where $V_1=\{x_{1},\ldots,x_{n}\}$ and $V_{2}=\{y_{1},\dots,y_{m}\}$. Therefore,
$I_{t}(G)=I_{q}J_{q+1}+I_{q+1}J_{q}$ if $t=2q+1$ and $I_{t}=I_{q}J_{q}$ if $t=2q$.
Thus Theorem \ref{qr} together with Theorem \ref{polymatroidal} now yields $I_{t}(G)$ is a ploymatroidal ideal.
}
\end{Example}

\medskip
Next we study the classes of generalized mixed product ideals that they are not polymatroidal.

\begin{Proposition}
\label{notpolymatroidal}
Let $L=I_{q}J_{r}+I_{s}J_{t} \subset K[x_1,\ldots, x_n,y_1, \ldots, y_m]$ be the generalized mixed product ideal, where for integers $a$ and $b$, the ideal $I_a$  (resp.\ $J_b$) is the ideal generated by all squarefree monomials of degree $a$ in the polynomial ring $K[x_1,\ldots, x_n]$
(resp.\ of degree $b$ in the polynomial ring $K[y_1, \ldots, y_m]$), and where $0 \leq s < q \leq n$, $0 \leq r < t \leq m$.
Then $L$ is not polymatroidal if $L$ can be written in one of the following forms:

{\em (a)} $L=I_{r}+J_{r}$, $r>1$.

{\em (b)} $L=J_{r}+I_{s}J_{t}$, with $2\leq s <n$ and $r=s+t$.

{\em (c)} $L=I_{q}J_{r}+I_{s}J_{t}$, with $s>q+1$, $q\geq1$, $t\geq1$ and $q+r=s+t$.
\end{Proposition}

\begin{proof}
(a) Let $u, v$ be two monomials belonging to $G(L)$, where $u=x_{i_{1}}\dots x_{i_{r}} \in G(I_{r})$ and $v=y_{j_{1}}\dots y_{j_{r}} \in G(J_{r})$ such that
$i_{1} < i_{2} < \dots <i_{r}$ and $j_{1}< j_{2}<\dots < j_{r}$, and where $\deg_{x_{i_{1}}}(u) >\deg(x_{i_{1}})(v)$.
Then for all $l$ with $1\leq l \leq r$ we have $\deg_{y_{j_{l}}}(u) <\deg y_{j_{l}} (v)$ and $y_{j_{l}}(u/x_{i_{1}})\notin G(L)$.
Therefore, $L$ is not polyamatroidal.

(b) Let $f=x_1\dots x_{s+1} y_{1}^{2}\dots y_{t-1}^{2}y_{t}\dots y_{r}$. We set $y_{0}=0$.
It then follows from \cite[Proposition 2.6]{RV} that $f\in \overline{L^{2}}$. By counting degrees, and using $s \geq 2$, we have $f\notin L^{2}$.
Therefore $L^{2}$ is not integrally closed, and hence $L$ is not polymatroidal by \cite[Theorm 3.4]{HRV}.

(c) The assertion follows by \cite[Proposition 2.7]{RV} and \cite[Theorm 3.4]{HRV}.
\end{proof}

\medskip
The {\em Veronese ideal} of degree $r$ is the ideal $I_{r}$ of $K[x_{1},\dots,x_{n}]$ which
is generated by all the monomials in the variables $x_{1},\ldots,x_{n}$ of degree
$r$: $I_{r}=(x_{1},\ldots,x_{n})^{r}$.

In the following, we consider the case that all $L_{ij}$ are powers of variables.

\begin{Theorem}
\label{expansion}
Let $L=I_{q}J_{r}+I_{q+1}J_{r-1}\subset K[x_1,\ldots, x_n,y_1, \ldots, y_m]$ be the generalized mixed product ideal, where for integers $a$ and $b$, the ideal $I_a$  (resp.\ $J_b$) is Veronese ideal of degree $a$ in the polynomial ring $K[x_1,\ldots, x_n]$
(resp.\ of degree $b$ in the polynomial ring $K[y_1, \ldots, y_m]$), and where $q \geq0$, $r\geq 1$. Then $L$ is polymatroidal.
\end{Theorem}

\begin{proof}
Let $u=x_{1}^{a_1}\dots x_{n}^{a_n} y_{1}^{b_1}\dots y_{m}^{b_m}$ and $v=x_{1}^{c_1}\dots x_{n}^{c_n} y_{1}^{d_1}\dots y_{m}^{d_m}$ be two monomials in $G(L)$. We set $a=\sum_{i=1}^{n}a_{i}$, $b=\sum_{j=1}^{m}b_{j}$, $c=\sum_{i=1}^{n}c_{i}$ and $d=\sum_{j=1}^{m}d_{j}$.
We may suppose that $a_{i}>c_{i}$.
If $a \leq c$, then there exists some $k$ such that $a_{k}< c_{k}$. Therefore, $x_{k}(u/x_{i})\in G(L)$ and the assertion follows. If $a>c$, thus there exists $j$ such that $b_{j}< d_{j}$. Hence, $y_{j}(u/x_{i})\in G(L)$ which is desired assertion.

Now, we assume that $b_{j} > d_{j}$. If $b\leq d$, then there exists some $l$ such that $b_{l}< d_{l}$, and $y_{l}(u/y_{j})\in G(L)$, as desired.
If $b> d$, then there exists some $i$ such that $a_{i} <c_{i}$. Thus $x_{i}(u/y_{j})\in G(L)$, and the conclusion follows.
\end{proof}

\begin{Theorem}
\label{two}
Let $L=I_{q}J_{r}\subset K[x_1,\ldots, x_n,y_1, \ldots, y_m]$ be the generalized mixed product ideal, where for integers $a$ and $b$, the ideal $I_a$  (resp.\ $J_b$) is Veronsese ideal generated by all monomials of degree $a$ in the polynomial ring $K[x_1,\ldots, x_n]$
(resp.\ of degree $b$ in the polynomial ring $K[y_1, \ldots, y_m]$).
Then $L$ is a polymatroidal ideal if and only if one of the following conditions are satisﬁed:

{\em (1)} $q=0$;

{\em (2)} $r=0$;

{\em (3)} $q,r> 0$.
\end{Theorem}

\begin{proof}
The cases $q=0$ or $r=0$ being clear, we assume that $q>0$ and $r>0$.
The ideal $L$ is a generalized mixed product ideal induced by the ideal $I=(x^{q}y^{r}).$
Let $u=x_{1}^{a_{1}}\dots x_{n}^{a_{n}} y_{1}^{b_{1}}\dots y_{m}^{b_{m}}$ and $v=x_{1}^{c_{1}}\dots x_{n}^{c_{n}} y_{1}^{d_{1}}\dots y_{m}^{d_{m}}$ be two monomials in $G(L)$. We may assume that $a_{i}> c_{i}$. Then there exists some $l$ such that $a_{l}< c_{l}$. Therefore, $x_{l}(u/x_{i})\in G(L)$.
Now suppose that $b_{j}> d_{j}$. Thus there exists some $k$ such that $b_{k}< d_{k}$ and $y_{k}(u/y_{j})\in G(L)$, as desired.
\end{proof}

\begin{Proposition}
\label{notnormal}
Let $L=I_{q}J_{r}+I_{s}J_{t}\subset K[x_1,\ldots, x_n,y_1, \ldots, y_m]$ be the generalized mixed product ideal, where for integers $a$ and $b$, the ideal $I_a$  (resp.\ $J_b$) is Veronese ideal of degree $a$ in the polynomial ring $K[x_1,\ldots, x_n]$
(resp.\ of degree $b$ in the polynomial ring $K[y_1, \ldots, y_m]$), and where $0 \leq s < q \leq n$, $0 \leq r < t \leq m$.
Then $L$ is not polymatroidal if $L$ can be written in one of the following forms:

{\em (a)} $L=I_{r}+J_{r}$, $r>1$.

{\em (b)} $L=J_{r}+I_{s}J_{t}$, with $2\leq s < n$ and $r=s+t$.
\end{Proposition}

\begin{proof}
{\em (a)} The ideal $L$ is a generalized mixed product ideal induced by $I=(x^{r},y^{r})$ with $r> 1$.
Let $u,v$ be two monomials in $G(L)$, where $u=x_{1}^{a_{1}}\dots x_{n}^{a_{n}} \in G(I_{r})$ and $v=y_{1}^{b_{1}}\dots y_{n}^{b_{n}} \in G(J_{r})$, and where $\deg_{x_{i}}(u) >\deg_{x_{i}}(v)$.
Therefore for all $j$ with $1\leq j \leq r$ we have $\deg_{y_{j}}(u) <\deg y_{j} (v)$ and $y_{j}(u/x_{i})\notin G(L)$.
Then $L$ is not polyamatroidal.

{\em (b)} Let $u=x_{1}^{s}y_{1}^{t}$ and $v=y_{1}^{r}$ be two monomials belonging to $G(L)$, where $\deg_{x_{1}}(u)>\deg_{x_{1}}(v)$. Therefore, $\deg_{y_{1}}(u)<\deg_{y_{1}}(v)$ such $y_{1}(x_{1}^{s}y_{1}^{t}/x_{1})\notin G(L)$, and hence $L$ is not polymatroidal.
\end{proof}

\medskip

Let $I$ be an ideal in a ring $R$. An element $r\in R$ is said to be {\em integral over I} if there exists
an integer $n$ and elements $a_{i}\in I^{i}, i=1,\dots,n$ such that
\[
r^{n}+a_{1}r^{n-1}+a_{2}r^{n-2}+\dots+a_{n-1}r=0
\]
Such an equation is called {\em an equation of integral dependence of r over I (of degree n)}.

The {\em integral closure} of $I$ is the set of all elements of $R$ which are integral over $I$.
The integral closure of a monomial ideal is again a monomial ideal.

The set of all elements that are integral over $I$ is called {\em the integral closure} of $I$,
and is denoted $\overline{I}$. In general $I\subseteq \overline{I}$.
If $I=\overline{I}$, $I$ is said to be {\em integrally closed} or {\em complete}.
If all the powers $I^{k}$ are integrally closed, $I$ is said to be {\em normal}.

In \cite{MN} we observed how the generalized mixed product ideal commutes with the integral closure of a monomial ideal.
Furthermore, we gave a geometric description of the integral closure of $L$.

\begin{Theorem}
\label{normality}
\cite[Theorem 2.4]{MN}
Let $L=\sum_{j=1}^m\prod_{i=1}^{n}L_{i,\ab_j(i)}$ be a generalized mixed product ideal,
induced by the monomial ideal $I$ with $G(I)=\{\xb^{\ab_1},\ldots,\xb^{\ab_m}\}$,
where the ideals $L_{i,\ab_j(i)}$ are Veronese ideals of degree $\ab_j(i)$ in the variables $x_{i1},x_{i2},\ldots,x_{im_i}$.
Then $I$ is normal if and only if $L$ is normal.
\end{Theorem}

\begin{Example}
\label{principal}
{\em
Let $L=L_{1,a_1}L_{2,b_1}$ with $a_1,b_1>1$
be the generalized mixed product ideal
induced by a monomial ideal $I=(x_1^{a_1}x_2^{b_1})$, where
the ideal  $L_{1,a_1}$  (resp.\ $L_{2,b_1}$) is the ideal generated by all monomials of degree $a_1$ in the polynomial ring $ K[x_{11},\ldots, x_{1m_1}]$
(resp.\ of degree $b_1$ in the polynomial ring $ K[x_{21}, \ldots, x_{2m_2}]$).
We have
\[
\overline{(x_1^{a_1})}\cap \overline{(x_2^{b_1})}=(x_1^{a_1})\cap (x_2^{b_1})\subseteq\overline{(x_1^{a_1})\cap (x_2^{b_1})}.
\]

For all $r\in \overline{(x_1^{a_1})\cap (x_2^{b_1})}$ there exist an equation
$r^{n}+c_{1} r^{n-1}+ \dots +c_{n-1} r + c_n=0$ with $c_i \in ((x_1^{a_1})\cap (x_2^{b_1}))^{i}$ for $i=1,\ldots,n$.
It follows that $c_i\in (x_1^{a_1})^i$ and $c_i\in  (x_2^{b_1})^i$.
Thus $r\in \overline{(x_1^{a_1})}\cap \overline{(x_2^{b_1})}$.
Since $\gcd((x_1^{a_1}), (x_2^{b_1}))=1$, it then follows that
$(x_1^{a_1}) (x_2^{b_1})$ is a $\lcm ((x_1^{a_1}), (x_2^{b_1}))$,
thus $(x_1^{a_1}) (x_2^{b_1})\in (x_1^{a_1})\cap (x_2^{b_1})$. Then $I$
is complete.
For all $k>0$, we have
\[
I^k=(x_1^{a_1})^k (x_2^{b_1})^k=(x_1^{a_1})^k\cap (x_2^{b_1})^k,
\]
hence $I^k$ is integrally closed.
Therefore Theorem ~\ref{normality} implies that
\[
L=L((x_1^{a_1} x_2^{b_1}); \{L_{ij}\})=L_{1,a_1}L_{2,b_1}
\]
is normal.
}
\end{Example}

In the following, we study the normality of generalized mixed product ideals induced by a  monomial ideal in $K[x_1,x_2]$,
where $L_{ij}$ are squarefree Veronese ideals in $K[x_{i1},x_{i2},\ldots,x_{im_i}]$ of degree $a_{j}$.

\begin{Theorem}
\label{normality}
Let $L=\sum_{r=1}^sI_{a_r}J_{b_r}\subset K[x_1,\ldots, x_n,y_1, \ldots, y_m]$ where the ideals $I_{a_r}$ in
$K[x_{1},\ldots,x_{n}]$ and the ideals $J_{b_r}$ in
$K[y_{1},\ldots,y_{m}]$ are squarefree Veronese ideals of degree $a_r$ and $b_r$,
respectively. Assume further that $0\leq a_1< \dots <a_s\leq n$, $m\geq b_1>\dots> b_s\geq0$ and $a_{r}+b_{r}=h$ for all $r=1,\dots s$.
If $a_{r+1}=a_{r}+1$ and $b_{r-1}=b_{r}+1$ for all $r=1,\dots,s-1$ then $L$ is normal.
\end{Theorem}

\begin{proof}
The ideal $L$ is a generalized mixed product ideal induced by the ideal
$$I=(x^{a_1}y^{b_1},\dots ,x^{a_s}y^{b_s}).$$
Let $u=x_{1}^{a_{c1}}\dots x_{n}^{a_{cn}} y_{1}^{b_{c1}}\dots y_{m}^{b_{cm}}$ and $v=x_{1}^{a_{d1}}\dots x_{n}^{a_{dn}} y_{1}^{b_{d1}}\dots y_{m}^{b_{dm}}$ be two monomials in $G(L)$ for some $c,d$ where $c,d=1,\dots,r$.
We set $a_{r}=\sum_{i=1}^{n}a_{ri}$ and $b_{r}=\sum_{j=1}^{m}b_{rj}$ for all $r=1,\dots,s$.

Suppose that $a_{ci}> a_{di}$. If $a_{c}\leq a_{d}$, then there exists some $l$ such that $a_{cl} < a_{dl}$ and $x_{l}(u/x_{i})\in G(L)$.
If $a_{c} > a_{d}$, then $b_{c} < b_{d}$. Furthermore, $b_{c} < b_{d}$ implies that there is $k$ with $b_{ck} < b_{dk}$. Therefore $y_{j}(u/x_{i})\in G(L)$.

We may assume that $b_{cj} > b_{dj}$. If $b_{c} \leq b_{d}$, then there exists some $w$ such that $b_{cw} < b_{dw}$. Hence, $y_{w}(u/y_{j})\in G(L)$.
If $b_{c} > b_{d}$, then there exists some $i$ such that $a_{ci}< a_{di}$. Therefore it is clear that $x_{i}(u/y_{j})\in G(L)$. Thus we conclude that $L$ is a polymatroidal ideal. Hence, by applying \cite[Theorem 3.4]{HRV} we obtain that $L$ is normal. This yields the desired conclusion.
\end{proof}

\section{Linear quotients of generalized mixed product ideals}
\label{three}
Let $K$ be a field and $S=K[x_1,\dots,x_n]$ the polynomial
over $K$ in the variables $x_1,\ldots,x_n$, and let $I\subset S$ be a monomial ideal with $I\neq S$ whose minimal set of generators  is $G(I)=\{\xb^{\ab_1},\ldots, \xb^{\ab_m}\}$.

Let $K$ be a field and $T=K[x_{11},\ldots,x_{1m_1},x_{21},\ldots,x_{2m_2},\ldots,x_{n1},\ldots,x_{nm_n}]$
be the polynomial ring over $K$ in the variables
\[
x_{11},\ldots,x_{1m_1},x_{21},\ldots,x_{2m_2},\ldots,x_{n1},\ldots,x_{nm_n},
\]
and $L$ be defined as in $(\ref{lj})$.
The main goal of this section is to study the linear quotients of generalized mixed product ideals.

A monomial ideal $I$ has {\em linear quotients} if the monomials that minimally generate $I$ can be ordered $g_{1},\dots, g_{q}$ such that for all
$1\leq i\leq q-1$, $((g_{1},\dots,g_{i}):g_{i+1})$ is generated by linear forms $x_{i_{1}},\dots,x_{i_{t}}$.

Now we consider the generalized mixed product ideal $L$ induced by a monomial ideal $I$, where the ideals $L_{ij}$ are squarefree Veronese ideals in $K[x_{i1},x_{i2},\ldots,x_{im_i}]$ of degree $a_{j}$.

\begin{Theorem}
\label{linear}
Let $L=\sum_{r=1}^sI_{a_r}J_{b_r}\subset K[x_1,\dots,x_n,y_1,\dots,y_m]$ be the generalized mixed product ideal induced by a monomial ideal
$I=(x^{a_1}y^{b_1},\dots,x^{a_s}y^{b_s})$, where the ideals $I_{a_r}$ in $K[x_{1},\ldots,x_{n}]$ and the ideals $J_{b_r}$ in $K[y_{1},\ldots,y_{m}]$ are squarefree Veronese ideals of degree $a_r$ and $b_r$, respectively.
Assume further that $a_s> \dots >a_1 :=1$, $b_1>\dots> b_s :=1$ and $m,n> 1$. If

\item [(1)] $a_{r+1}=a_{r}+1$ and $b_{r-1}=b_{r}+1$ for all $r=1,\dots,s-1$,

\item [(2)] $a_{r}+b_{r}=h$ for all $r=1,\dots s$ with $2 \leq h \leq m+n-1$

then $L$ have linear quotients.
\end{Theorem}

\begin{proof}
Let $u\in G(L)$.
The ideal $L$ is a generalized mixed product ideal induced by the ideal
\[
I=(x^{a_1}y^{b_1},\dots ,x^{a_s}y^{b_s}),
\]
where the ideals $I_{a_r}$ in
$K[x_{1},\ldots,x_{n}]$ and the ideals $J_{b_r}$ in
$K[y_{1},\ldots,y_{m}]$ are squarefree Veronese ideals of degree $a_r$ and $b_r$,
respectively. Suppose that $m,n> 1$.
We set $L^{*}=(v \in G(L): v\prec u)$ with $\prec$ the lexicographical order on $x_{1},\ldots,x_{n},y_{1},\dots,y_{m}$ induced by
\[
x_{1}\succ x_{2}\succ \cdots \succ x_{n}\succ y_{1}\succ y_{2}\succ \cdots\succ y_{m}.
\]
We claim that $L^{*}:u=(v/\gcd(u,v):v\in L^{*})$ is generated by variables. Thus we must prove that for any $v\prec u$, there exists a variable of
$K[x_{1},\ldots,x_{n},y_{1},\dots,y_{m}]$ in $L^{*}:u$ such that it divides $v/\gcd(u,v)$.

Let $u=x_{1}^{a_1}\dots x_{n}^{a_n} y_{1}^{b_1}\dots y_{m}^{b_m}$ and $v=x_{1}^{c_1}\dots x_{n}^{c_n} y_{1}^{d_1}\dots y_{m}^{d_m}$ be two monomials in $G(L)$. We know that $v \prec u$, there exists an integer $i$ with $a_{i}> c_{i}$ and $a_{l}=c_{l}$ for $l=1,\dots ,i-1$. Therefore, there exists an integer $j$
with $a_{j}< c_{j}$ such that $x_{j}(u/x_{i})\in G(L)$. Since $j< i$ it follows that $x_{j}(u/x_{i})\in L^{*}$. Then $x_{j}\in L^{*}:u$.
Since the $j$-th component of the vector exponent of $v/\gcd(u,v)$ is given by $c_{j}-\min \{c_{j},a_{j}\}> 0$, we conclude that $x_{j}$
divides $v/\gcd(u,v)$ as desired.

If $a_{l}=c_{l}$ for all $l=1,\dots,n$, $b_{i}> d_{i}$ and $b_{t}=d_{t}$ for all $t=1,\dots, i-1$, $i\in \{1,\dots ,m\}$, hence $y_{j}\in L^{*}:u$
and $y_{j}$ divides $v/\gcd(u,v)$. This yields the
desired conclusion.
\end{proof}

\begin{Example}
\label{L3}
{\em
Let $L=L_{1,1}L_{2,1}+L_{1,2}L_{2,1}\subset T=K[x_{11},\ldots, x_{1m_{1}},x_{21}, \ldots, x_{2m_{2}}]$ be the generalized mixed product ideal
induced by a monomial ideal $I=(xy^2,x^2y)$, where for integers $a$ and $b$, the ideal  $L_{1,a}$  (resp.\ $L_{2,b}$) is the ideal generated by all squarefree monomials of degree $a$ in the polynomial ring $K[x_{11},\ldots, x_{1m_{1}}]$
(resp.\ of degree $b$ in the polynomial ring $K[x_{21}, \ldots, x_{2m_{2}}]$), and where $m_{1},m_{2} > 1$.
Therefore $L$ have linear quotients, by Theorem \ref{linear}.
}
\end{Example}

In the following, we consider the case that the generalized mixed product ideal
induced by a Veronese ideal.

\begin{Theorem}
\label{veronese}
Let $L$ be a generalized mixed product ideal,
induced by the monomial ideal $I$ with $G(I)=\{\xb^{\ab_1},\ldots,\xb^{\ab_m}\}$, where the ideals $L_{i,\ab_j(i)}$ are Veronese ideals of degree $\ab_j(i)$
in the variables $x_{i1},x_{i2},\ldots,x_{im_i}$.
If $I$ is a Veronese ideal then $L$ have linear quotients.
\end{Theorem}

\begin{proof}
For a monomial ideal $I\subset S$ with $G(I)=\{\xb^{\ab_1},\ldots,\xb^{\ab_m}\}$, let
\[
L(I;\{L_{ij}\})=\sum_{j=1}^m\prod_{i=1}^{n}L_{i,\ab_j(i)}.
\]
It follows from \cite[Theorem 2.3]{M} that $L(I;\{L_{ij}\})^k$ is a generalized mixed product ideal induced by $I^k$ and $L(I;\{L_{ij}\})^k=L(I^{k};\{L_{ij}\})$
for all $k\geq1$.
Let $I $ be a Veronese ideal of degree $k$.
Then $I=(x_{1},\ldots,x_{n})^{k}$, which
is generated by all the monomials in the variables $x_{1},\ldots,x_{n}$ of degree $k$.
Then by \cite[Theorem 2.3]{M} we have
\begin{eqnarray*}
L(I;\{L_{ij}\})
&=&L((x_{1},\ldots,x_{n})^{k};\{L_{ij}\})\\
&=&  L((x_{1},\ldots,x_{n});\{L_{ij}\})^{k}\\
&=& (x_{11},\ldots,x_{1m_1},x_{21},\ldots,x_{2m_2},\ldots,x_{n1},\ldots,x_{nm_n})^{k}.
\end {eqnarray*}

Hence \cite[Corollary 12.6.4]{HH} implies that $L$ have linear quotients.
Thus the desired conclusion follows.
\end{proof}

\begin{Theorem}
\label{mixed}
Let $L=\sum_{j=1}^mL_{1,a_j}L_{2,b_j}$ be the generalized mixed product ideal
induced by a monomial ideal
$I=(x_1^{a_1}x_2^{b_1},\dots ,x_1^{a_m}x_2^{b_m}),$
where the ideals $L_{1,a_j}$ in
$K[x_{11},x_{12},\ldots,x_{1m_1}]$ and the ideals $L_{2,b_j}$ in
$K[x_{21},x_{22},\ldots,x_{2m_2}]$ are polymatroidal ideals of degree $a_j$ and $b_j$,
respectively. Assume that $0\leq a_1< \dots <a_m\leq m_1$ and $m_2\geq b_1>\dots> b_m\geq0$.
Furthermore, let $I$ has a linear resolution.
Then $L$ has a linear resolution.
\end{Theorem}

\begin{proof}
The assertion follows by \cite[Theorem 2.3]{HMY}.
\end{proof}

\medskip

Next assume in addition that $I\subset S$ is generated in one degree and suppose that $I$ has linear quotients with respect to the ordering $g_{1},\dots,g_{q}$ of the monomials belonging to $G(I)$. Then the colon ideal $(g_{1},g_{2},\dots ,g_{i}):g_{i+1}$ is generated by a subset of $\{x_{1},x_{2},\dots,x_{n}\}$
for each $1\leq i\leq q-1$. Let $r_{i}$ denote the number of variables which is required to generate $(g_{1},g_{2},\dots g_{i}):g_{i+1}$.
Let $r(I)=\max_{1 \leq i\leq q-1} r_{i}$.

\begin{Theorem}
\label{r(I)}
Let $T=K[x_{11},\ldots,x_{1m_1},x_{21},\ldots,x_{2m_2}]$ with $m_{1},m_{2}>1$.
Let $L=\sum_{j=1}^mL_{1,a_j}L_{2,b_j}\subset T$ where the ideals $L_{1,a_j}$ in
$K[x_{11},x_{12},\ldots,x_{1m_1}]$ and the ideals $L_{2,b_j}$ in
$K[x_{21},x_{22},\ldots,x_{2m_2}]$ are squarefree Veronese ideals of degree $a_j$ and $b_j$,
respectively. Assume further that $a_m> \dots >a_1 :=1$, $b_1>\dots> b_m :=1$ and $m> 1$. If

\item[(1)] $a_{j+1}=a_{j}+1$ and $b_{j-1}=b_{j}+1$ for all $j=1,\dots,m-1$,

\item[(2)] $a_{j}+b_{j}=z$ for all $j=1,\dots m$ with $2 < z \leq m_{1}+m_{2}-1$

then
\[
r(\sum_{j=1}^mL_{1,a_j}L_{2,b_j})=m_{1}+m_{2}-1.
\]
\end{Theorem}

\begin{proof}
Let $L=\sum_{j=1}^mL_{1,a_j}L_{2,b_j} \subset T$ be the generalized mixed product ideal induced by a monomial ideal
\[
I=(x^{a_1}y^{b_1},\dots,x^{a_m}y^{b_m}),
\]
where the ideals $L_{1,a_j}$ in $K[x_{11},x_{12},\ldots,x_{1m_1}]$ and the ideals $L_{2,b_j}$ in
$K[x_{21},x_{22},\ldots,x_{2m_2}]$ are squarefree Veronese ideals of degree $a_j$ and $b_j$,
respectively.

We order the generators of $L$ with respect to the monomial order $\prec_{lex}$ on the variables $x_{11},x_{12},\dots,x_{1m_{1}},x_{21},x_{22},\dots,x_{2m_{2}}$
induced by
\[
x_{11}\succ x_{12}\succ \cdots\succ x_{1m_{1}}\succ x_{21} \succ x_{22}\succ \cdots \succ x_{2m_{2}}.
\]
Now we assume that $l\in \{1,\dots,m_{1}\}$ such that $(l-1)+a_{1l}=z-1$, this implies that
\[
x_{11}x_{12}\cdots x_{1l-1}x_{1l}^{a_{1l}}x_{21}\succ x_{11}x_{12}x_{1l-1}x_{1l}^{a_{1l}}x_{22}\succ \cdots \succ x_{11}x_{12}x_{1l-1}x_{1l}^{a_{1l}}x_{2m_{2}}\succ \cdots
\]
and so on up to $x_{1m_{1}}x_{2k}^{a_{2k}}x_{2k+1}\dots x_{2m_{2}}$ with $(m_{2}-k)+a_{2k}=z-1$.
It then follows that the maximum system of their generators is $\{x_{11},\ldots,x_{1m_1},x_{21},\ldots,x_{2m_2-1}\}$.
Therefore, $r(L)=m_{1}+m_{2}-1$.
\end{proof}

\medskip
There are two important invariants attached to a graded ideal $I \subset S$
defined in terms of the minimal graded free resolution of $S/I$.

We consider the minimal free graded resolution of $M=S/I$ as an $S$-module.
\[
0\rightarrow  \oplus_{j\in \mathbf{z} } S(-j)^{\beta_{pj}}(M)\rightarrow \dots \rightarrow \oplus_{j\in \mathbf{Z}  } S(-j)^{\beta_{0j}}(M)\rightarrow M\rightarrow 0
\]
The {\em Castelnuovo-Mumford regularity} (or simply the regularity) of $M =S/I$ is defined as
\[
\reg(S/I) := max\{j - i : \beta_{i,j} \neq 0\}.
\]
Furthermore, the {\em projective dimension} of $M$ is defined as
\[
\pd(M) := max\{i :\beta_{i,j} \neq 0  \quad  \text{for some} \quad j\}.
\]

Now we investigate algebraic invariants of $L$.

\begin{Theorem}
\label{pd}
Let $T=K[x_{11},\ldots,x_{1m_1},x_{21},\ldots,x_{2m_2}]$ with $m_{1},m_{2}>1$.
Let $L=\sum_{j=1}^mL_{1,a_j}L_{2,b_j}\subset T$ where the ideals $L_{1,a_j}$ in
$K[x_{11},x_{12},\ldots,x_{1m_1}]$ and the ideals $L_{2,b_j}$ in
$K[x_{21},x_{22},\ldots,x_{2m_2}]$ are squarefree Veronese ideals of degree $a_j$ and $b_j$,
respectively. Assume further that $a_m> \dots >a_1 :=1$, $b_1>\dots> b_m :=1$ and $m> 1$. If

\item[(1)] $a_{j+1}=a_{j}+1$ and $b_{j-1}=b_{j}+1$ for all $j=1,\dots,m-1$,

\item[(2)] $a_{j}+b_{j}=z$ for all $j=1,\dots m$ with $2 < z \leq m_{1}+m_{2}-1$

then
\[
\pd(T/L)=m_{1}+m_{2}.
\]
\end{Theorem}

\begin{proof}
According to \cite[Corollary 1.6]{HT} we have the length of the minimal free resolution of $T/L$ over $T$ is equal to $r(L)+1$.
Hence by applying Theorem \ref{r(I)} we obtains that $\pd(T/L)=m_{1}+m_{2}-1+1=m_{1}+m_{2}.$
\end{proof}

\medskip

Recall that if a monomial ideal generated in the same degree has linear quotients, hence it has a linear resolution, (see \cite{CH}).

In the following, we compute the regularity of $L$ in the case that all $L_{ij}$ are squarefree Veronese ideals
in $K[x_{i1},x_{i2},\ldots,x_{im_i}]$ of degree $a_{j}$.

\begin{Theorem}
\label{reg}
Let $T=K[x_{11},\ldots,x_{1m_1},x_{21},\ldots,x_{2m_2}]$ with $m_{1},m_{2}>1$.
Let $L=\sum_{j=1}^mL_{1,a_j}L_{2,b_j} \subset T$ be the generalized mixed product ideal induced by a monomial ideal
\[
I=(x^{a_1}y^{b_1},\dots,x^{a_m}y^{b_m}),
\]
where the ideals $L_{1,a_j}$ in $K[x_{11},x_{12},\ldots,x_{1m_1}]$ and the ideals $L_{2,b_j}$ in
$K[x_{21},x_{22},\ldots,x_{2m_2}]$ are squarefree Veronese ideals of degree $a_j$ and $b_j$,
respectively. Assume further that $a_m> \dots >a_1 :=1$, $b_1>\dots> b_m :=1$ and $m> 1$. If

\item[(1)] $a_{j+1}=a_{j}+1$ and $b_{j-1}=b_{j}+1$ for all $j=1,\dots,m-1$,

\item[(2)] $a_{j}+b_{j}=z$ for all $j=1,\dots m$ with $2 < z \leq m_{1}+m_{2}-1$

then
\[
\reg(T/L)=z-1.
\]
\end{Theorem}

\begin{proof}
By Theorem \ref{linear}, we have $L$ has a linear quotients. Therefore, $L$ has a linear resolution and hence $\reg(L)=z$, as desired.
\end{proof}

\begin{Example}
\label{L4}
{
For a monomial ideal $I=(xy^{3},x^{2}y^{2},x^{3}y)$ let $L=L_{1,1}L_{2,3}+L_{1,2}L_{2,2}+L_{1,3}L_{2,1}\subset T=K[x_{11},x_{12}, x_{13},x_{21},x_{22}, x_{23}]$ be the generalized mixed product ideal, where for integers $a$ and $b$, the ideal $L_{1,a}$  (resp.\ $L_{2,b}$) is the ideal generated by all squarefree monomials of degree $a$ in the polynomial ring $K[x_{11},x_{12}, x_{13}]$
(resp.\ of degree $b$ in the polynomial ring $K[x_{21},x_{22}, x_{23}]$).
Therefore,
\begin{eqnarray*}
L&=&(x_{11}x_{21}x_{22}x_{23}, x_{12}x_{21}x_{22}x_{23}, x_{13}x_{21}x_{22}x_{23},x_{11}x_{12}x_{21}x_{22}, x_{11}x_{12}x_{21}x_{23}, \\
&&x_{11}x_{12}x_{22}x_{23},x_{11}x_{13}x_{21}x_{22},x_{11}x_{13}x_{21}x_{23}, x_{11}x_{13}x_{22}x_{23}, x_{12}x_{13}x_{21}x_{22},\\
&& x_{12}x_{13}x_{21}x_{23},x_{12}x_{13}x_{22}x_{23}, x_{11}x_{12}x_{13}x_{21}, x_{11}x_{12}x_{13}x_{22},x_{11}x_{12}x_{13}x_{23}).
\end{eqnarray*}
Hence, by applying Theorems \ref{reg} and \ref{pd} we obtain that $reg(T/L)=3$ and  $\pd(T/L)=6$.
}
\end{Example}

\end{document}